\def\twobytwo#1#2#3#4{\pmatrix{ #1&#2\cr #3&#4\cr}}
\centerline{\bf Eigenfunctions on the Finite Poincar\'e Plane}
\centerline{\it by \rm Jinghua Kuang
\footnote{${}^*$}{Research at MSRI supported in part by NSF grant \#DMS 9022140}}
\bigskip
\bf 1. Introduction. \rm Let $F_q$ be a finite field with $q$ elements ($q$ odd.) Fix a non-square elemnet $\delta\in F$. 
$$
H_q=F_q(\sqrt{\delta})-F_q=\{x+y\sqrt{\delta}\,\vert\,x,y\ne0\in F_q\}
$$
is called the finite Poincar\'e plane. From [1], [2], [4], [7], we have the following facts:

\item{1.} $G={\rm GL}_2(F_q)$ acts on $H_q$  by linear transformation:
$$
gz={{az+b}\over{cz+d}},\qquad\forall g=\twobytwo abcd\in G,\,\forall z\in H_q.
$$
Let $K=\left\{\twobytwo a{b\delta}ba\right\}$. One may identify $G/K$ with $H_q$ and with $P=\left\{\twobytwo yx01\big\vert\,y\ne 0, x\in F_q\right\}$.

\item{2.} There is a $G$-invariant pseudo-distance $\Delta$ on $H_q$:
$$
\Delta(z_1,\,z_2)={{N(z_1-z_2)}\over{{\rm Im}z_1\cdot{\rm Im}z_2}}\qquad\forall z_1,z_2\in H_q,
$$
where $N$ is the norm from $F_q(\sqrt{\delta})$ to $F_q$ and ${\rm Im}z=y$ if
$z=x+y\sqrt{\delta}$.

\item{3.} For $a\in F^\times_q$ and $a\ne 4\delta$, a graph structure 
$X_q(\delta, a)$ is defined on $H_q$: $(z_1, z_2)$ is an edge iff $\Delta(z_1,z_2)=a$. Katz [3] and Li [5] proved that $X_q(\delta, a)$ are Ramanujan graphs using different methods. 

However, distributions of eigenvalues of the adjacency matrices have not been well understood. We will calculate the first and second moments of the asymptotic distribution of the eigenvalues of the adjacency matrices and provide evidence for Terras' conjecture in [7].

Let $V$ be the space of all complex functions on $H_q$. Using character sums, Evans [2] constructed an orthogonal basis of $V$ that diagonalizes all the adjacency matrices. Kuang [4], using group representation theory, also constructed such an orthogonal basis. This paper will also compare these two bases. 

\smallskip
\bf 2. Hecke Algebra. \rm Since $K$ is the analogue of O(2) in the classical Poincar\'e upper half plane, it is natural to consider the Hecke algebra $H(G, K)$,
which is defined as the algebra of all bi-$K$-invariant complex functions on $G$ under convolution. Since $K$ is the isotropic subgroup of $\sqrt{\delta}$,
$\Delta(g,1)=\Delta(g^{-1},1)$ implies that $KgK=Kg^{-1}K$. Therefore, $H(G, K)$ is commutative. Obviously, $H(G, K)$ has dimension $q$. Hence there exist $q$ idempotents in $H(G, K)$.

\smallskip
\bf 3. Idempotents. \rm For an irreducible representation $\pi$ of $G$, let $\eta_\pi'(g)={1\over{|K|}}\sum_{k\in K}{\rm tr\,}\pi(kg)$. $\pi$ is spherical (i.e. ${\rm Res}^G_K(\pi)\supset 1_K$) iff $\eta_\pi'(g)\ne 0$. From [4], we have
$$
\eta_{\pi_1}'*\eta_{\pi_2}'=
\left\{\matrix{
0,\hfill  & {\rm if } \pi_1\not\simeq\pi_2;\hfill\cr
{  {|G|}  \over {{\rm dim}(\pi)}  }\eta_\pi', &{\rm if }\pi_1\simeq\pi_2\simeq\pi.
}\right. 
$$
 Let $\omega$ be a multiplicative character of $F_q(\sqrt{\delta})$ of order $q^2-1$, $\chi$ be a multiplicative character of order $q-1$, and $s=\chi^{(q-1)/2}$. There are $q$ irreducible spherical representations of $G$ (see [6]):
$$
\eqalign{
\pi&=1_G;\cr
{\rm or }&=\rho(s,s)=q\,{\rm dimensional\, irreducible\, component\, of\, }{\rm Ind}^G_B(s,s);\cr
{\rm or }&=\rho(\chi^j,\bar{\chi}^j)={\rm Ind}^G_B(\chi^j,\bar{\chi}^j),\quad j=1,\dots,{{q-3}\over 2};\cr
{\rm or }&=\rho_{\nu_j}={\rm cuspidal\, representation\, corresponding\, to\, }\nu_j=\omega^{j(q-1)},\quad j=1,\dots,{{q-1}\over2}.\cr}
$$
($B$ is the Borel subgroup of $G$.) Let $\eta_\pi={{{\rm dim}(\pi)}\over{|G|}}\eta_\pi'$. Let $\eta_0=\eta_{1_G}$, $\eta_j=\eta_{\rho(\chi^j,\bar{\chi}^j)},\,j=1,\dots,(q-3)/2$, $\eta_{{q-1}\over2}=\eta_{\rho(s,s)}$ and $\eta_{{{q-1}\over2}+j}=\eta_{\rho_{\nu_j}},\,j=1,\dots,(q-1)/2.$ Then

\proclaim Theorem 1 (Kuang [4]). $\eta_0,\eta_1,\dots,\eta_{q-1}$ are the $q$ idempotents in $H(G, K)$.

\smallskip
\bf 4. Representation of $H(G, K)$. \rm $H(G, K)$ has a natural representation on $L^2(H_q)=L^2(P)$, the space of all complex functions on $H_q=P$: 
$$
T_\varphi(f)(p)={1\over{|K|}}f*\varphi(p)
={1\over{|K|}}\sum_{h\in G}f(h)\varphi(h^{-1}p)
=\sum_{p_1\in P}f(p_1)\varphi(p_1^{-1}p),
$$ 
for $\varphi\in H(G, K)$, $f\in L^2(P)$ and $p\in P$.

For the usual inner product $\langle,\,\rangle$ on $L^2(P)$, we have
$\langle T_\varphi(f_1),f_2\rangle=\langle f_1,T_{\hat{\varphi}}(f_2)\rangle$,
where $\hat{\varphi}(h)=\overline{\varphi(h^{-1})}$. That is, $H(G, K)$ acts on $L^2(P)$ self-adjointly. Hence, there exists a basis of simultaneous eigenfunctions of $H(G, K)$.
\smallskip
\bf 5. Construction of the Bases. \rm Let $\psi$ be a fixed non-trivial additive character of $F_q$. Denote $\psi_a(x)=\psi(ax)$. Define
\item{}$\chi_i(g)=\chi^i(y)$ if $gK=\twobytwo yx01 K$, $i=1,\dots,q-1$.
\item{}$\Psi_a(g)=\delta(y)\psi_a(x)$ if $gK=\twobytwo yx01 K$, $a\in F_q^\times$,

\noindent where $\delta(y)=1$ if $y=1$; $=0$ otherwise. Then
\proclaim Theorem 2 (Kuang [4]). $\chi_i,\, i=1,\dots,q-1$ and $\Psi_a*\eta_i,\,i=1,\dots,q-1,\,a\in F_q^\times$ are simultaneous eigenfunctions of $H(G, K)$ and make up an orthogonal basis of $L^2(P)$. Moreover,
$T_{\eta_j}(\chi_i)=\chi_i$ if $\eta_j=\eta_{\rho(\chi^{\pm i},\bar{\chi}^{\pm i})}$; $=0$ otherwise. $T_{\eta_j}(\Psi_a*\eta_i)=\Psi_a*\eta_i$ if $i=j$; $=0$ otherwise.

Terras [8] constructed the following functions:
$$
K_{i,\psi_a}(\twobytwo yx01)=\psi_a(-x)\chi^i(y)\sum_{u\in F_q}\psi_a(u)\bar{\chi}^i(u^2-\delta y^2).
$$
Evans [2], based on Velasquez's work according to [9], constructed the following functions:
$$
H_{t, i, \psi_a}(\twobytwo yx01)={\psi_a(-x)\over q+1}\sum_{u\in F_q}{}^*
\psi_a(u)\sum_{\alpha\in U_1}\omega^i(\alpha)
s\left(\alpha+{1\over\alpha}+{\delta y\over t}
+{t\over\delta y}-{u^2\over ty}\right),
$$
where $U_1=\{\alpha\in F_q^\times\vert\, N(\alpha)=1\}$ and the asterisk indicates that when $y=\pm t$, the $u=0$ term is multiplied by $q+1$.

\proclaim Theorem 3 (Terras, Velasquez, Evans [2]). $K_{i,\psi_a}$, $i=1,\dots, q-1$, $a\in F_q$ and $H_{t,i,\psi_a}$, $t\in F_q^\times$, $i=1,\dots, q-1$, $a\in F_q$
are simultaneous eigenfunctions of $H(G, K)$. There exist a $t_0$ for any fixed $i$ and $a$ such that $K_{i,\psi_0}$, $i=1,\dots, q-1$, 
$K_{i,\psi_a}$, $i=1,\dots,(q-1)/2$, $a\in F_q$ and $H_{t_0, i, \psi_a}$ $i=1,\dots,(q-1)/2$, $a\in F_q^\times$ make up an orthogonal basis of $L^2(P)$.

\smallskip
\bf 6. Comparison of the Bases. \rm We now prove that the two bases constructed above are essentially same. Let us denote by $\pi_i$ the spherical representation that give rise to $\eta_i$. 

\proclaim Theorem 4. The following identities hold. 
$$
\eqalign{
\overline{C_i(0)}\cdot\chi_i&=K_{i, 1},\, i=1,\dots,q-1;\cr
\Psi_a*\eta_i&={{\rm dim}(\pi_i)\over |G|(q+1)}C_i(a)K_{i,\psi_{-a}},\,i=1,\dots,(q-1)/2;\cr
\Psi_a*\eta_{i+(q-1)/2}&={{\rm dim}(\pi_{i+(q-1)/2})\over |G|}H_{-\delta,i,\psi_{-a}},\,i=1,\dots,(q-1)/2,\cr
}
$$
where $C_i(a)=\sum_{x\in F_q}\psi_a(x)\chi^i(x^2-\delta)$.

\noindent \bf Proof. \rm $K_{i,1}=\overline{C_i(0)}\cdot\chi_i$ is obvious. Since for any pair of $a, b\in F_q$, $K_{i,\psi_b}$ and $\Psi_a*\eta_i$ are eigenfunctions of $H(G, K)$ of the same eigenvalues ($i=1,\dots,(q-1)/2$), $\Psi_a*\eta_i$ is a complex linear combination of $K_{i,\psi_b},\,(b\in F_q)$. Now, 
$$
\eqalign{
\langle K_{i,\psi_b},\Psi_a*\eta_i\rangle&=\sum_{p\in P}K_{i,\psi_b}(p)
\overline{\Psi_a*\eta_i}(p)\cr
&=\sum_{p\in P}\sum_{u\in F_q}\psi_b(u-x)\chi^i(y)\overline{\chi}^i(u^2-\delta y^2)
\overline{\eta_i}(p)\sum_{v\in F_q}\psi(-(a+b)v)\cr
&=0\quad{\rm if }\,a+b\ne 0\cr
}
$$
Hence $\Psi_a*\eta_i$ is a multiple of $K_{i,\psi_{-a}}$. Using the equation (2.17) in [2], we calculate
$$
\eqalign{
\Psi_a*\eta_i(1)&=\sum_{v\in F_q}\psi_a(u)\eta_i(\twobytwo 1{-u}01)\cr
&={{\rm dim}(\pi_i)\over |G|(q+1)}C_i(a)\overline{C_i(a)}.\cr
}
$$
And $K_{i, \psi_{-a}}(1)=\overline{C_i(a)}$.  Therefore,
$\Psi_a*\eta_i={{\rm dim}(\pi_i)\over |G|(q+1)}C_i(a)K_{i,\psi_{-a}}$.

For the third identity, we use the equation (2.16) in [2]. Let $|S(z)|$ be the cardinality of the $K$-orbit of $z$ in $H_q$. Let $\hat{i}=i+(q-1)/2$. Then,
$$
\eqalign{
\Psi_a*\eta_{\hat{i}}(p)&=\sum_{u\in F_q}\psi_a(u)\eta_{\hat{i}}(\twobytwo y{x-u}01)\cr
&=\sum_{u\in F_q}\psi_a(x-u)\eta_{\hat{i}}(\twobytwo yu01)\cr
&=\psi_a(x)\sum_{u\in F_q}\psi_a(-u){{\rm dim}(\pi_{\hat{i}})\over |G|}{1\over |S(u+y\sqrt{\delta})|}\sum_{\alpha\in U_1}\omega^i(\alpha)s\!\left(
\alpha+{1\over \alpha}-2+{\Delta(u+y\sqrt{\delta},\sqrt{\delta})\over \delta}\right)\cr
&=\psi_a(x)\sum_{u\in F_q}\psi_a(-u){{\rm dim}(\pi_{\hat{i}})\over |G|}{1\over |S(u+y\sqrt{\delta})|}\sum_{\alpha\in U_1}\omega^i(\alpha)s\!\left(
\alpha+{1\over \alpha}+{u^2\over y\delta}-y-{1\over y}\right)\cr
&={{\rm dim}(\pi_{\hat{i}})\over |G|}H_{-\delta,i,\psi_{-a}}(p)\cr
}
$$
\smallskip

\bf 7. Comments and Questions. \rm For any function $f$ on $G$, define the Fourier coefficient
$$
F(g;a,f)=\sum_{x\in F_q}f(\twobytwo 1x01 g)\psi_a(x).
$$
From [4], we see
$\langle \Psi_a*\eta_i,\Psi_a*\eta_i\rangle=|K|\cdot q\cdot F(1;a,\eta_i).$
So $F(1;a,\eta_i)\ne 0$. But can one evaluate $F(1;a,\eta_i)$? We have $F(1;a,\eta_i)=\Psi_{-a}*\eta_i(1)$. Hence, for $i=1,\dots,(q-1)/2$,
$$
\eqalign{
F(1;a,\eta_i)&={{\rm dim}(\pi_i)\over |G|(q+1)}|C_i(-a)|^2
={{\rm dim}(\pi_i)\over |G|(q+1)}|C_i(a)|^2\cr
F(1;a,\eta_{\hat{i}})&={{\rm dim}(\pi_{\hat{i}})\over |G|}\sum_{u\in F_q}\psi_a(u)
{1\over |S(u+\sqrt{\delta})|}\sum_{\alpha\in U_1}\omega^i(\alpha)s\!\left(\alpha+
{1\over\alpha}+{u^2\over \delta}-2\right)\cr
}
$$
However, we don't know much about $C_i(a)$ and  about the last set of sums. We note that $C_i(0)$ is essentially a Jacobi sum, in fact, 
$C_i(0)=-\chi^i(-\delta)J(\chi^i,s)$.
\smallskip

\bf 8. Connection to $X_q(\delta,a)$. \rm There are $q$ double cosets: $D_a=KgK$, where $\Delta(g,1)=a\in F_q$. Let $\varphi_a$ be the characteristic function of the set $D_a$. Then $A_a={1\over |K|}T_{\varphi_a}$ acting on $L^2(H_q)$ is the adjacency matrix of $X_q(\delta,a)$. Therefore each basis constructed in Section 5 diagonalizes all the adjacency matrices. Let $S_a=D_a\cap P$. Then
$$
\lambda_i(a)={|G|\over{\rm dim}(\pi_i)}|S_a|\eta_i(D_a),\,i=0,1,\dots,q-1
$$
are all the eigenvalues of $A_a$.
Katz [3], Li [5] proved that $|\lambda_i(a)|\leq 2\sqrt{q}$ for $i=1,\dots,q-1$, which confirms that $X_a(\delta,a)$ is Ramanujan.

Fix $a\ne 4\delta\in F^\times$, Terras [7] conjectured that 
$\{\lambda_i(a)/\sqrt q\,\vert\,i=1,\dots,q-1\}$ asymptotically has Sato-Tate distribution, i.e. for $E\subset [-2,2]$,
$$
\lim_{q\to\infty}{1\over q-1}|\{\lambda_i(a)\,\vert\,\lambda_i(a)/\sqrt q\in E\}|
={1\over 2\pi}\int_E\sqrt{4-x^2}dx.
$$
\smallskip

\bf 9. Moments. \rm Define a $q\times q$ matrix
$$
M=\left(\sqrt{|S_a|\over{\rm dim}(\pi_i)}\eta_i(S_a)\right)_{i=0,1,\dots,q-1,\,a\in F_q}.
$$
The idempotent property of $\eta_i$ ($i=0,\dots,q-1$) implies that $M M'={1\over |K||G|}I_q$. ($M'$ is the transpose of $M$.) Hence, $M'M={1\over|K||G|}I_q$, that gives
$$
\eqalign{
{1\over q-1}\sum_{i=0}^{q-1}{{\rm dim}(\pi_i)\over q}
\cdot{\lambda_i(a)\over\sqrt q}&=0\cr
{1\over q-1}\sum_{i=0}^{q-1}{{\rm dim}(\pi_i)\over q}
\cdot\left({\lambda_i(a)\over\sqrt q}\right)^2&={|G||S_a|\over |K|q^2(q-1)}\cr
}
$$
That, in turn, implies
$$
\eqalign{
\lim_{q\to\infty}{1\over q-1}\sum_{i=1}^{q-1}
{\lambda_i(a)\over\sqrt q}&=0\cr
\lim_{q\to\infty}{1\over q-1}\sum_{i=1}^{q-1}
\left({\lambda_i(a)\over\sqrt q}\right)^2&=1.\cr
}
$$
Therefore, the first and second moments of $\{\lambda_i(a)/\sqrt q\,\vert\,i=1,\dots,q-1\}$ asymptotically match with those of the Sato-Tate distribution.
\smallskip

\bf Acknowledgement. \rm This is a written-up of a talk that the author delivered at the Workshop on Spectral Theory of Automorphic Forms and Number Theory, in MSRI, Berkeley, California, October 24-28, 1994. Some calculation was completed while the author was visiting MSRI. He thanks the MSRI for its hospitality.
\bigskip

\bf References.\rm\hfill\break
\noindent [1]. \bf J. Angel, N. Celniker, S. Poulos, A. Terras, C. Trimble, E. Velasquez \it Special Functions on Finite Upper Half Planes, \rm in AMS Contemporary Math. Vol. 138, 1992.

\noindent [2]. \bf R. Evans \it Character sums as orthogonal eigenfunctions of adjacency operators for Cayley graphs, \rm in Proc. of a conference on finite fields, AMS contemporary Math. Vol. 168, 1994.

\noindent [3]. \bf N. Katz \it Estimates for Soto-Andrade Sums, \rm J. f\"ur die Reine und Angew. Math. 438 (1993), 143-161.

\noindent [4]. \bf J. Kuang \it A Natural Orthogonal Basis of Eigenfunctions of the Hecke Algebra Acting on Cayley Graphs, \rm to appear in the Proceeding of American Mathematical Society.

\noindent [5]. \bf W. Li \it A survey of Ramanujan graphs, \rm Proc. of a Conference on Arithmetic Geometry and Coding Theory at Luminy, 1993 
(to appear).

\noindent [6]. \bf I. Piatetski-Shapiro \it Complex Representation of ${\rm GL}_2(F_q)$, \rm AMS Contemporary Math. Vol. 16, 1983.

\noindent [7]. \bf A. Terras \it Are Finite Upper Half Plane Graphs Ramanujan? \rm in DIMACS Series in Discrete Mathematics and Theoretic Computer Science, Vol. 10 (1993), 125-142

\noindent [8]. \bf A. Terras \it Fourier Analysis on Finite Groups and Applications, \rm U. C. S. D. Lecture Notes, 1991-2.

\noindent [9]. \bf E. Velasquez, \rm Private Conversation, October 27, 1994.

\bigskip
\bigskip
\noindent Department of Mathematics

\noindent Pennsylvania State University, at Fayette

\noindent Uniontown, PA 15401

\noindent E-mail address: \it kuang@math.psu.edu \rm

\bye